\documentclass[12pt]{article}
\usepackage{amssymb}

\usepackage{latexsym}
\usepackage{pictex}
\usepackage[mathscr]{eucal}
\usepackage{graphicx,color,psfrag}

\newtheorem{theorem}{Theorem}[section]

\renewcommand{\epsilon}{\varepsilon}

\begin{document}

\title{On the discrete spectrum of a spatial quantum waveguide with a disc window}
\author{H. Najar \thanks { D\'epartement de Math\'ematiques, ISMAI. Kairouan, Bd Assed Ibn Elfourat, 3100 Kairouan Tunisia.}
\and S. Ben Hariz  \thanks{\ D\'epartement de Math\'ematiques, Universit\'e du Maine, Le Mans, France.}
\and M. Ben Salah
 \thanks{\ D\'epartement de Math\'ematiques, ISMAI. Kairouan, Bd Assed Ibn Elfourat, 3100 Kairouan Tunisia.}}
\date{}
\maketitle

\begin{abstract}
In this study we investigate the bound states of the Hamiltonian
describing a quantum particle living on three dimensional straight
strip of width $d$. We impose the Neumann boundary condition on a
disc window of radius $a$ and Dirichlet boundary conditions on the
remained part of the boundary of the strip. We prove that such
system exhibits discrete eigenvalues below the essential spectrum
for any $a>0$. We give also a numeric estimation of the number of
discrete eigenvalue as a function of $\displaystyle \frac{a}{d}$.
When $a$ tends to the infinity, the asymptotic of the eigenvalue is
given.
\end{abstract}

\baselineskip=20pt \setcounter{section}{0}
\renewcommand{\theequation}{\arabic{section}.\arabic{equation}}
%%%%%%%%%%%%%%%%%%%%%%%%%%%%%%%%%%%%%%%%%%%%%%%%%%%%%%%%%%%%%%%%%

\noindent


\section{Itroduction}

The study of quantum waves on quantum waveguide has gained much interest and
has been intensively studied during the last years for their important
physical consequences. The main reason is that they represent an interesting
physical effect with important applications in nanophysical devices, but
also in flat electromagnetic waveguide. See the monograph \cite{hurt} and
the references therein. \newline
Exner et al. have done seminal works in this field. They obtained results in
different contexts, we quote \cite{exner1, exner4, exner2, exner3}. Also in
\cite{naj6,speis1,Stokel} research has been conducted in this area; the
first is about the discrete case and the two others for deals with the
random quantum waveguide.

It should be noticed that the spectral properties essentially
depends on the geometry of the waveguide, in particular, the
existence of a bound states induced by curvature
\cite{bulla,exner4,1,exner2} or by coupling of straight waveguides
through windows \cite{exner2, hurt} were shown. The waveguide with
Neumann boundary condition were also investigated in several papers
\cite{david, naz}. A possible next generalization are waveguides
with combined Dirichlet and Neumann boundary conditions on different
parts of the boundary. The presence of different boundary conditions
also gives rise to nontrivial properties like the existence of bound
states.The rest of the paper is organized as follows, in Section 2,
we define the model and recall some known results. In section 3, we
present the main result of this note followed by a discussion.
Section 4 is devoted for numerical experiments.

\section{The model}

The system we are going to study is given in Fig 1. We consider a
Schr\"odinger particle whose motion is confined to a pair of parallel plans
of width $d$. For simplicity, we assume that they are placed at $z=0$ and $%
z=d$. We shall denote this configuration space by $\Omega$
\[
\Omega=\Bbb{R}^2\times [0,d].
\]
Let $\gamma(a)$ be a disc of radius $a$, without loss of generality
we assume that the center of $\gamma(a)$ is the point $(0,0,0)$;
\begin{equation}
\gamma(a)=\{(x,y,0)\in \Bbb{R}^3;\ x^2+y^2\leq a^2\}.
\end{equation}
We set $\Gamma=\partial \Omega\diagdown \gamma(a)$. We consider Dirichlet
boundary condition on $\Gamma$ and Neumann boundary condition in $\gamma(a)$.

\begin{figure}[h]
\centering
\resizebox{0.95\textwidth}{!}{\includegraphics{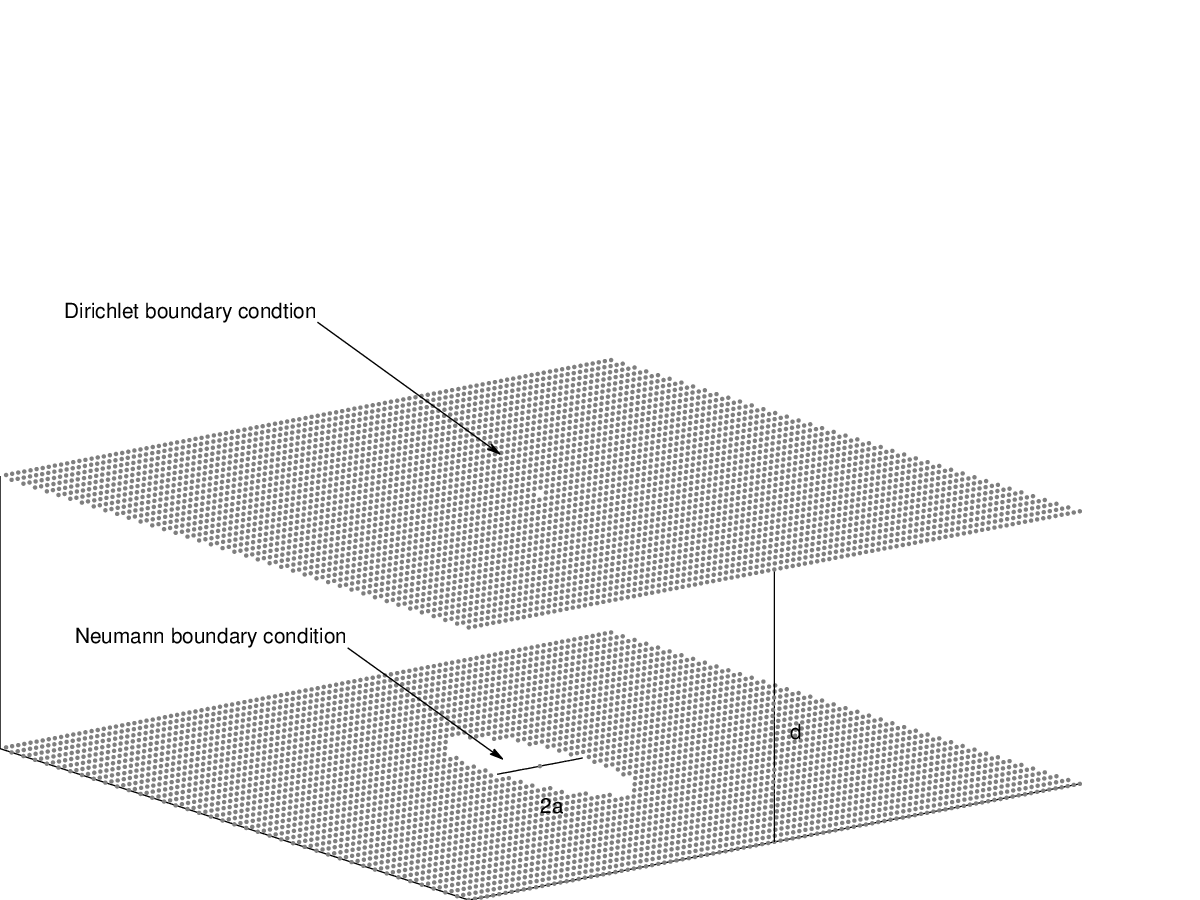}}
\caption{The waveguide with a disc window and  two different
boundaries conditions} \label{fig1}
\end{figure}

\subsection{The Hamiltonian}
Let us define the self-adjoint operator on $L^{2}(\Omega )$
corresponding to the particle Hamiltonian $H$. This is will be done
by the mean of quadratic forms. Precisely, let $q_0$ be the
quadratic form
\begin{equation}
q_{0}(f,g)=\int_{\Omega }\nabla f\cdot \overline{\nabla g}d^{3}x,\ \mathrm{%
with\ domain}\ \mathcal{Q}(q_{0})=\{f\in H^{1}(\Omega );\ f\lceil \Gamma
=0\},
\end{equation}
where $H^{1}(\Omega )=\{f\in L^{2}{(\Omega )}|\nabla f\in L^{2}(\Omega )\}$
is the standard Sobolev space and we denote by $f\lceil \Gamma $, the trace
of the function $f$ on $\Gamma $. It follows that $q_{0}$ is a densely
defined, symmetric, positive and closed quadratic form. We denote the unique
self-adjoint operator associated to $q_{0}$ by $H$ and its domain by $%
D(\Omega )$. It is the hamiltonian describing our system. From \cite{Resi}
(page 276), we infer that the domain $D(\Omega )$ of $H$ is
\[
D(\Omega )=\Big\{f\in H^{1}(\Omega );\ -\Delta f\in L^{2}(\Omega ),f\lceil
\Gamma =0,\frac{\partial f}{\partial z}\lceil \gamma (a)=0\Big\}
\]
and
\[
Hf=-\Delta f,\ \ \forall f\in D(\Omega ).
\]

\subsection{Some known facts}

Let us start this subsection by recalling that in the particular case when $%
a=0$, we get $H^{0}$, the Dirichlet Laplacian, and $a=+\infty $ we get $%
H^{\infty }$, the Dirichlet-Neumann Laplacian. Since
\[
H=(-\Delta _{\Bbb{R}^{2}})\otimes I\oplus I\otimes (-\Delta _{\lbrack 0,d]}),%
\mathrm{on}\ L^{2}(\Bbb{R}^{2})\otimes L^{2}([0,d]),
\]
( see \cite{Resi}) we get that the spectrum of $H^{0}$ is $[(\frac{\pi }{2d}%
)^{2},+\infty \lbrack $. Consequently, we have
\[
\left[ (\frac{\pi }{d})^{2},+\infty \right[ \subset \sigma (H)\subset \left[
(\frac{\pi }{2d})^{2},+\infty \right[ .
\]
Using the property that the essential spectra is preserved under compact
perturbation, we deduce that the essential spectrum of $H$  is

\[
\sigma _{ess}(H)=\left[ (\frac{\pi }{d})^{2},+\infty \right] .
\]
An immediate consequence is  the discrete spectrum lies in $\left[ (\frac{%
\pi }{2d})^{2},(\frac{\pi }{d})^{2}\right] $.

\subsection{Preliminary: Cylindrical coordinates}

Let us notice that the system has a cylindrical symmetry, therefore,  it is
natural to consider the cylindrical coordinates system $(r,\theta ,z)$.
Indeed, we have that

\[
L^{2}(\Omega ,dxdydz)=L^{2}(]0,+\infty \lbrack \times \lbrack 0,2\pi \lbrack
\times \lbrack 0,d],rdrd\theta dz),
\]
We note by $\dot{\langle}\dot{,}\rangle _{r}$, the scaler product in $%
L^{2}(\Omega ,dxdydz)=L^{2}(]0,+\infty \lbrack \times \lbrack 0,2\pi \lbrack
\times \lbrack 0,d],rdrd\theta dz)$ given by
\[
\langle f,g\rangle _{r}=\int_{]0,+\infty \lbrack \times \lbrack 0,2\pi
\lbrack \times \lbrack 0,d]}fgrdrd\theta dz.
\]
We denote the gradient in cylindrical coordinates by $\nabla _{r}$. While
the Laplacian operator in cylindrical coordinates is given by
\begin{equation}
\Delta _{r,\theta ,z}=\frac{1}{r}\frac{\partial }{\partial r}(r\frac{%
\partial }{\partial r})+\frac{1}{r^{2}}\frac{\partial ^{2}}{\partial \theta
^{2}}+\frac{d^{2}}{dz^{2}}.
\end{equation}
Therefore, the eigenvalue equation is given by
\begin{equation}
-\Delta _{r,\theta ,z}f(r,\theta ,z)=Ef(r,\theta ,z).  \label{vp}
\end{equation}
Since the operator is positive, we set $E=k^{2}$. The equation (\ref{vp}) is
solved by separating variables and considering $f(r,\theta ,z)=\varphi
(r)\cdot \psi (\theta )\chi (z).$ Plugging the last expression in equation (%
\ref{vp}) and first separate $\chi $ by putting all the $z$ dependence in
one term so that $\frac{\chi ^{\prime \prime }}{\chi }$ can only be
constant. The constant is taken as $-s^{2}$ for convenience. Second, we
separate the term $\frac{\psi "}{\psi }$ which has all the $\theta $
dependance. Using the fact that the problem has an axial symmetry and the
solution has to be $2\pi $ periodic and single value in $\theta $, we obtain
$\frac{\psi "}{\psi }$ should be a constant $-n^{2}$ for $n\in \Bbb{Z}$.
Finally, we get the following equation for $\varphi $
\begin{equation}
\varphi ^{\prime \prime }(r)+\frac{1}{r}\varphi ^{\prime }(r)+[k^{2}-s^{2}-%
\frac{n^{2}}{r^{2}}]\varphi (r)=0.  \label{bes}
\end{equation}
We notice that the equation (\ref{bes}), is the Bessel equation and \ its
solutions could be expressed in terms of Bessel functions. More explicit
solutions could be given by considering boundary conditions.

\section{The result}

The main result of this note is the following Theorem.

\begin{theorem}
\label{th1} The operator $H$ has at least one isolated eigenvalue in $\left[
(\frac{\pi }{2d})^{2},(\frac{\pi }{d})^{2}\right] $ for any $a>0$.

Moreover for $a$ big enough, if $\lambda (a)$ is an eigenvalue of
$H$ less then $\displaystyle \frac{\pi ^{2}}{d^{2}}$, then we have.
\begin{equation}
\lambda (a)=\left( \frac{\pi }{2d}\right) ^{2}+o\left( \frac{1}{a^{2}}%
\right) .  \label{as1}
\end{equation}
\end{theorem}

\textbf{Proof.} Let us start by proving the first claim of the
Theorem. To do so, we define the quadratic form $\mathcal{Q}_0$,

\begin{equation}
\mathcal{Q}_{0}(f,g)=\langle \nabla f,\nabla g\rangle _{r}=\int_{]0,+\infty
\lbrack \times \lbrack 0,2\pi \lbrack \times \lbrack 0,d]}(\partial _{r}f%
\overline{\partial _{r}g}+\frac{1}{r^{2}}\partial _{\theta }f\overline{%
\partial _{\theta }g}+\partial _{z}f\overline{\partial _{z}g})rdrd\theta dz,
\end{equation}
with domain
\[
\mathcal{D}_{0}(\Omega )=\left\{ f\in L^{2}(\Omega ,rdrd\theta dz);\nabla
_{r}f\in L^{2}(\Omega ,rdrd\theta dz);f\lceil \Gamma =0\right\} .
\]
Consider the functional $q$ defined by
\begin{equation}
q[\Phi ]=\mathcal{Q}_{0}[\Phi ]-(\frac{\pi }{d})^{2}\Vert \Phi \Vert
_{L^{2}(\Omega ,rdrd\theta dz)}^{2}.
\end{equation}
Since the essential spectrum of $H$ starts at $\displaystyle
(\frac{\pi }{d})^{2}$, if we
construct a trial function $\Phi \in \mathcal{D}_{0}(\Omega )$ such that $%
q[\Phi ]$ has a negative value then the task is achieved. Using the
quadratic form domain, $\Phi $ must be continuous inside $\Omega $ but not
necessarily smooth. Let $\chi $ be the first transverse mode, i.e.
\begin{equation}
\chi (z)=\left\{
\begin{array}{ccc}
\sqrt{\frac{2}{d}}\sin (\frac{\pi }{d}z) & \ \mathrm{if} & z\in (0,d) \\
0 & \mathrm{otherwise}. &
\end{array}
\right.
\end{equation}
For $\Phi (r,\theta ,z)=\varphi (r)\chi (z)$, we compute
\begin{eqnarray*}
q[\Phi ] &=&\langle \nabla _{r}\varphi \chi ,\nabla _{r}\varphi \chi \rangle
-(\frac{\pi }{d})^{2}\Vert \varphi \chi \Vert _{L^{2}(\Omega ,rdrd\theta
dz)}^{2}, \\
&=&\int_{]0,+\infty \lbrack \times \lbrack 0,2\pi \lbrack \times \lbrack
0,d]}\left( |\chi (z)|^{2}|\varphi ^{\prime }(r)|^{2}+|\varphi (r)||\chi
^{\prime }(z)|^{2}\right) rdrd\theta dz-(\frac{\pi }{d})^{2}\Vert \varphi
\chi \Vert _{L^{2}(\Omega ,rdrd\theta )}^{2} \\
&=&2\pi \Vert \varphi ^{\prime }\Vert _{L^{2}([0,+\infty \lbrack ,rdr)}^{2}
\end{eqnarray*}
Now let us consider an interval $J=[0,b]$ for a positive $b>a$ and a
function $\varphi \in \mathcal{S}\left( [0,+\infty \lbrack \right) $ such
that $\varphi (r)=1$ for $r\in J.$ We also  define a family $\{\varphi
_{\tau }:\tau >0\}$ by
\begin{equation}
\varphi _{\tau }(r)=\left\{
\begin{array}{ccc}
\varphi (r) & \ \mathrm{if} & r\in (0,b) \\
\varphi (b+\tau (\ln r-\ln b)) & \mathrm{if} & r\geq b.
\end{array}
\right.
\end{equation}
Let us write
\begin{eqnarray}
\Vert \varphi _{\tau }^{\prime }\Vert _{L^{2}([0,+\infty ),rdr)}
&=&\int_{(0,\infty )}|\varphi _{\tau }^{\prime }(r)|^{2}rdr ,  \nonumber \\
&=&\int_{(b,+\infty )}\tau ^{2}|\varphi ^{\prime }(b+\tau (\ln r-\ln
b))|^{2}rdr ,  \nonumber \\
&=&\int_{(b,+\infty )}\frac{\tau ^{2}}{r^{2}}|\varphi ^{\prime }(b+\tau (\ln
r-\ln b))|^{2}rdr,  \nonumber \\
&=&\tau \int_{(b,+\infty )}\frac{\tau }{r}|\varphi ^{\prime }(b+\tau (\ln
r-\ln b))|^{2}dr,  \nonumber \\
&=&\tau \int_{(0,+\infty )}|\varphi ^{\prime }(s)|^{2}ds=\tau \Vert \varphi
^{\prime }\Vert _{L^{2}((0,+\infty ))}^{2}.  \label{etoi}
\end{eqnarray}
Let $j$ be a localization function from $C_{0}^{\infty }(0,a)$ and for $\tau
,\varepsilon >0$ we define
\begin{equation}
\Phi _{\tau ,\varepsilon }(r,z)=\varphi _{\tau }(r)[\chi (z)+\varepsilon
j(r)^{2}]=\varphi _{\tau }(r)\chi (z)+\varphi _{\tau }\varepsilon
j^{2}(r)=\Phi _{1,\tau ,\varepsilon }(r,z)+\Phi _{2,\tau ,\varepsilon }(r).
\end{equation}
\begin{eqnarray*}
q[\Phi ] &=&q[\Phi _{1,\tau ,\varepsilon }+\Phi _{2,\tau ,\varepsilon }] \\
&=&\mathcal{Q}_{0}[\Phi _{1,\tau ,\varepsilon }+\Phi _{2,\tau ,\varepsilon
}]-(\frac{\pi }{d})^{2}\Vert \Phi _{1,\tau ,\varepsilon }+\Phi _{2,\tau
,\varepsilon }\Vert _{L^{2}(\Omega ,rdrd\theta dz)}^{2}. \\
&=&\mathcal{Q}_{0}[\Phi _{1,\tau ,\varepsilon }]-(\frac{\pi }{d})^{2}\Vert
\Phi _{1,\tau ,\varepsilon }\Vert _{L^{2}(\Omega ,rdrd\theta dz)}^{2}+%
\mathcal{Q}_{0}[\Phi _{2,\tau ,\varepsilon }]-(\frac{\pi }{d})^{2}\Vert \Phi
_{2,\tau ,\varepsilon }\Vert _{L^{2}(\Omega ,rdrd\theta dz)}^{2} \\
&+&2\langle \nabla _{r}\Phi _{1,\tau ,\varepsilon },\nabla _{r}\Phi
_{2,\tau ,\varepsilon }\rangle _{r}-(\frac{\pi }{d})^{2}\langle \Phi
_{1,\tau ,\varepsilon },\Phi _{2,\tau ,\varepsilon }\rangle _{r}.
\end{eqnarray*}
Using the properties of $\chi $, noting that the supports of $\varphi $ and $%
j$ are disjoints and taking into account equation (\ref{etoi}), we get
\begin{equation}
q[\Phi ]=2\pi \tau \Vert \varphi ^{\prime }\Vert _{L^{2}(0,+\infty )}-8\pi
d\varepsilon \Vert j^{2}\Vert _{L^{2}(0,+\infty )}^{2}+2\varepsilon ^{2}\pi
\{2\Vert jj^{\prime }\Vert _{(L^{2}(0,\infty ),rdr)}^{2}-(\frac{\pi }{d}%
)^{2}\Vert j^{2}\Vert _{(L^{2}(0,\infty ),rdr)}^{2}\}.  \label{cle}
\end{equation}
Firstly, we notice that only the first term of the last equation depends on $%
\tau $. Secondly, the linear term in $\varepsilon $ is negative and
could be chosen sufficiently small so that it dominates over the
quadratic one. Fixing this $\varepsilon $ and then choosing $\tau $
sufficiently small the right hand side of (\ref{cle}) is negative.
This ends the proof of the first claim.\newline
The proof of the second claim is based on bracketing argument. Let us split $%
L^{2}(\Omega ,rdrd\theta dz)$ as follows, $L^{2}(\Omega ,rdrd\theta
dz)=L^{2}(\Omega _{a}^{-},rdrd\theta dz)\oplus L^{2}(\Omega
_{a}^{+},rdrd\theta dz)$, with
\begin{eqnarray*}
\Omega _{a}^{-} &=&\{(r,\theta ,z)\in \lbrack 0,a]\times \lbrack 0,2\pi
\lbrack \times \lbrack 0,d]\}, \\
\Omega _{a}^{+} &=&\Omega \backslash \Omega _{a}^{-}.
\end{eqnarray*}
Therefore
\[
H_{a}^{-,N}\oplus H_{a}^{+,N}\leq H\leq H_{a}^{-,D}\oplus H_{a}^{+,D}.
\]
Here we index by $D$ and $N$ depending on the boundary conditions considered
on the surface $r=a$. The min-max principle leads to

\[
\sigma _{ess}(H)=\sigma _{ess}(H_{a}^{+,N})=\sigma _{ess}(H_{r}^{+,D})=\left[
(\frac{\pi }{d})^{2},+\infty \right[ .
\]
Hence if $H_{r}^{-,D}$ exhibits a discrete spectrum below $\displaystyle \frac{\pi ^{2}}{%
d^{2}}$, then $H$ do as well. We mention that this is not a necessary
condition. If we denote by $\lambda _{j}(H_{a}^{-,D}),\lambda
_{j}(H_{a}^{-,N})$ and $\lambda _{j}(H)$, the $j$-th eigenvalue of $%
H_{a}^{-,D}$, $H_{a}^{-,N}$ and $H$ respectively then, again the minimax
principle yields the following
\begin{equation}
\lambda _{j}(H_{a}^{-,N})\leq \lambda _{j}(H)\leq \lambda _{j}(H_{a}^{-,D})
\label{es1}
\end{equation}
and for $2\geq j$
\begin{equation}
\lambda _{j-1}(H_{a}^{-,D})\leq \lambda _{j}(H)\leq \lambda
_{j}(H_{a}^{-,D}).
\end{equation}
\newline
$H_{a}^{-,D}$ has a sequence of eigenvalues \cite{abra,wat}, given
by
\[
\lambda _{k,n,l}=\left( \frac{(2k+1)\pi }{2d}\right) ^{2}+\left( \frac{%
x_{n,l}}{a}\right) ^{2}.
\]
Where $x_{n,l}$ is the $l$-th positive zero of Bessel function of
order $n$ ( see \cite{abra,wat}) . The condition
\begin{equation}
\lambda _{k,n,l}<\frac{\pi ^{2}}{d^{2}},\label{gaza}
\end{equation}
yields that $k=0$, so we get
\[
\lambda _{0,n,l}=\left( \frac{\pi }{2d}\right) ^{2}+\left( \frac{x_{n,l}}{a}%
\right) ^{2}.
\]
This yields that the condition (\ref{gaza}) to be fulfilled, will
depends on the value of $\displaystyle \left(
\frac{x_{n,l}}{a}\right) ^{2}$.\newline We recall that $x_{n,l}$ are
the positive zeros of the Bessel function $J_{n}
$. So, for any $\lambda (a)$, eigenvalue of $H$, there exists, $%
n,l,n^{\prime },l^{\prime }\in \Bbb{N}$, such that
\begin{equation}
\frac{\pi ^{2}}{4d^{2}}+\frac{x_{n,l}^{2}}{a^{2}}\leq \lambda (a)\leq \frac{%
\pi ^{2}}{4d^{2}}+\frac{x_{n^{\prime },l^{\prime }}^{2}}{a^{2}}.
\end{equation}
The proof of (\ref{as1}) is completed by observing by that $x_{n,l}$ and $%
x_{n^{\prime },l^{\prime }}$ are independent from $a$. In Figure 2, the
domain of existence of $\lambda _{1}(H),\lambda _{2}(H)$ and $\lambda _{3}(H)
$ are represented.$\boxdot $

\section{Numerical computations}

This section is devoted to some numerical computations. In \cite{glo} and
\cite{mar}, the number of positive zeros of Bessel functions less than $%
\lambda $ is estimated by $\displaystyle \frac{\lambda ^{2}}{\pi
^{2}}$ which is based on the approximate formula for the roots of
Bessel functions for large $l$ is
\begin{equation}
x_{n,l}\sim (n+2l-\frac{1}{2})\frac{\pi }{2}.
\end{equation}
taking into account $(\ref{es1})$, we get that for, $d$ and $a$
positives such that $\displaystyle \frac{a^{2}}{d^{2}}<\lambda
^{\ast }=1,9276$, $H$ has a unique discrete eigenvalue.

\begin{figure}[h]
\centering
\resizebox{0.99\textwidth}{!}{\includegraphics{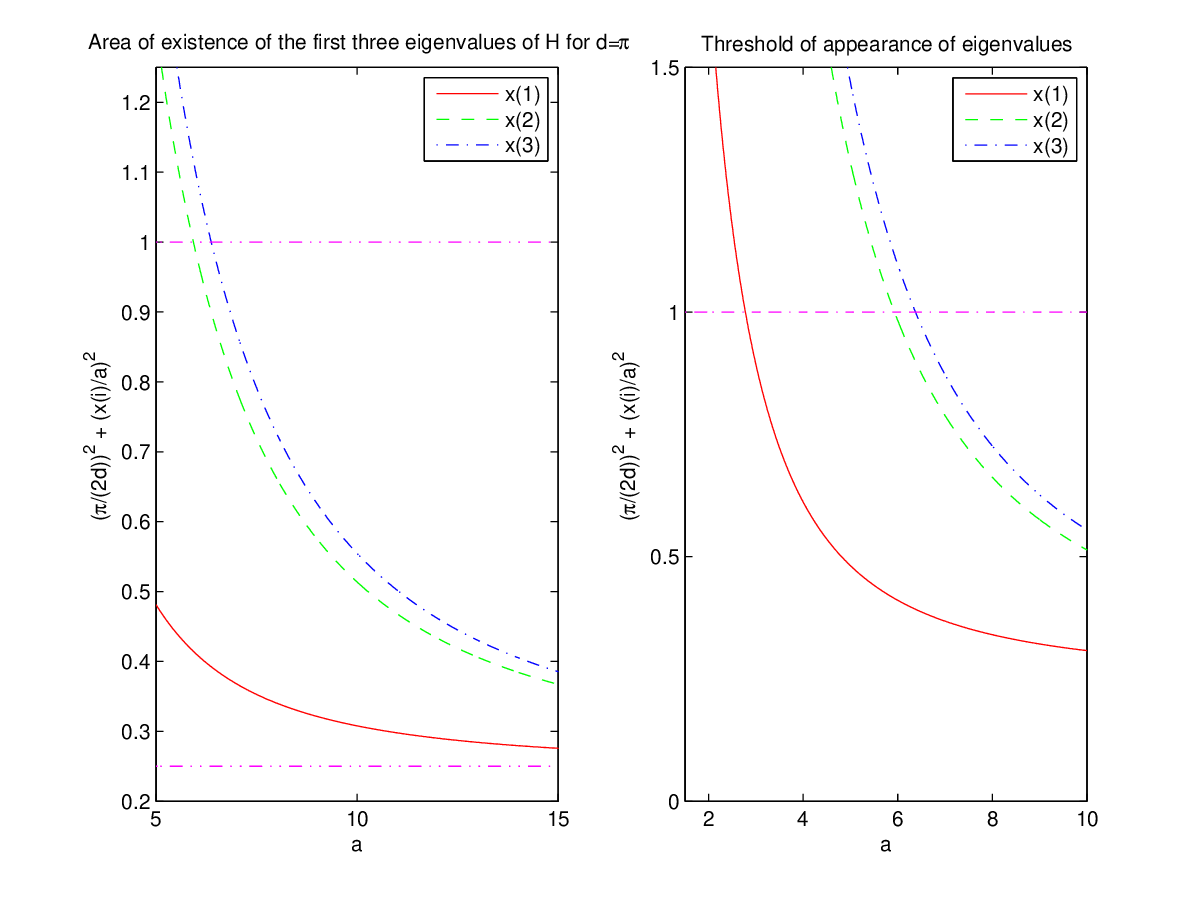}}
\caption{We represent $a \mapsto (\frac{\protect\pi}{2d})^2+(\frac{x(i)}{a}%
)^{2}$ where $x(1), x(2), x(3)$ are the first three zeros of the bessel
functions increasingly ordered. }
\label{fig2}
\end{figure}

\begin{figure}[h]
\centering
\resizebox{0.5\textheight}{!}{\includegraphics{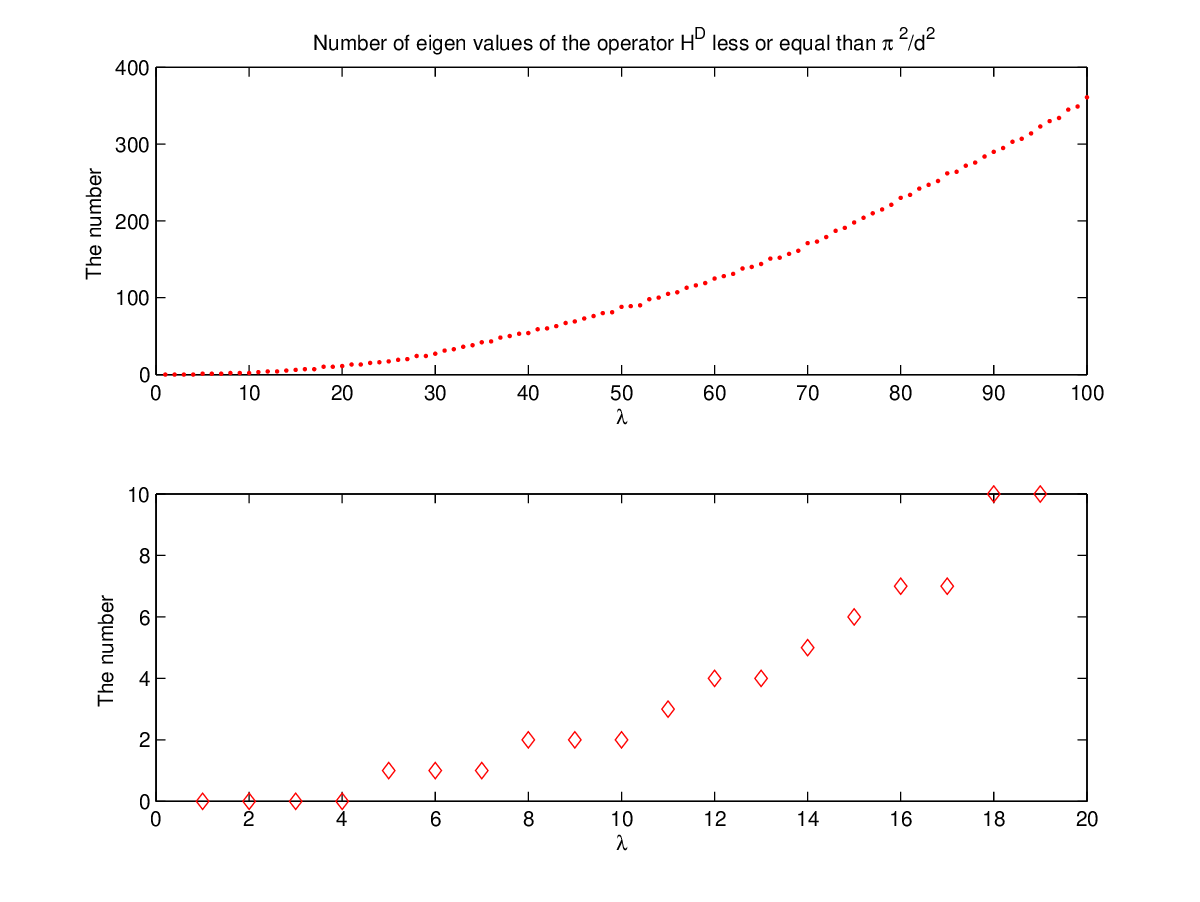}}
\caption{ The number of the eigenvalues of the operator $H^D$  function of $\lambda\equiv a/d.$
%\newline $\protect\lambda\mapsto
%\sum_{n=0}^{\protect\sqrt{\protect\pi}(2d)^{-1}
%\protect\lambda} \sum_{l\geq 1}I_{\{x_{n,l}<\protect\sqrt{\protect\pi}%
%(2d)^{-1} \protect\lambda\} }$
}
\label{fig3}
\end{figure}

\begin{figure}[h]
\centering
\resizebox{0.5\textheight}{!}{\includegraphics{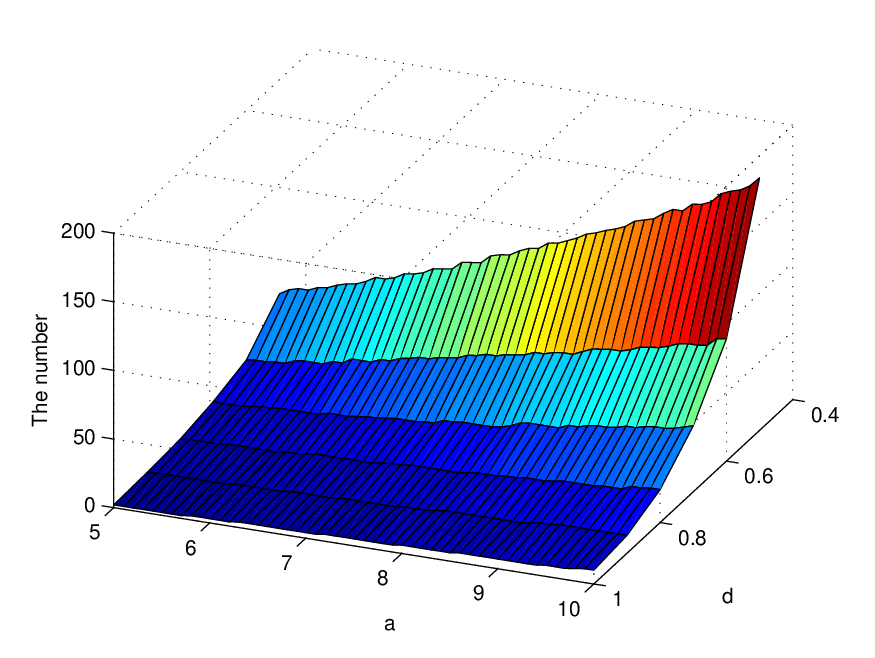}}
\caption{ The number of the eigenvalues of the operator $H^D$  function of $d$ and $a.$
%\newline $\protect\lambda\mapsto
%\sum_{n=0}^{\protect\sqrt{\protect\pi}(2d)^{-1}
%\protect\lambda} \sum_{l\geq 1}I_{\{x_{n,l}<\protect\sqrt{\protect\pi}%
%(2d)^{-1} \protect\lambda\} }$
}
\label{fig4}
\end{figure}

\clearpage

\newpage

\noindent\textbf{Acknowledgements.} It is a pleasure for the authors to
thanks Professor Pavel Exner for useful discussions, valuable comments and
remarks which significantly improve this work.

\end{document}